\documentclass[10pt,english]{article}
\usepackage[T1]{fontenc}
\usepackage[latin1]{inputenc}
\usepackage{a4wide}
\usepackage{babel}
\usepackage{latexsym}
\usepackage{amsmath}
\usepackage{amssymb}
\usepackage{hyperref}
\usepackage{color}

\makeatletter

\newtheorem{prop}{Proposition}[section]
\newtheorem{lemma}[prop]{Lemma}
\newtheorem{theorem}[prop]{Theorem}
\newtheorem{defi}[prop]{Definition}
\newtheorem{rem}[prop]{\em Remark\/}
\newtheorem{cor}[prop]{Corollary}
\newtheorem{cond}[prop]{Conditions}
\newtheorem{ex}[prop]{\em Example\/}
\newtheorem{sch}[prop]{Scholium}

\newcommand{\bp}{\begin{prop}}
\newcommand{\bl}{\begin{lemma}}
\newcommand{\bt}{\begin{theorem}}
\newcommand{\bd}{\begin{defi}\rm}
\newcommand{\br}{\begin{rem}\rm}
\newcommand{\be}{\begin{equation}}
\newcommand{\bea}{\begin{eqnarray}}
\newcommand{\bpr}{\begin{proof}}
\newcommand{\bc}{\begin{cor}}
\newcommand{\bco}{\begin{cond}}
\newcommand{\bex}{\begin{ex}\rm}
\newcommand{\bsch}{\begin{sch}}

\newcommand{\ep}{\end{prop}}
\newcommand{\el}{\end{lemma}}
\newcommand{\et}{\end{theorem}}
\newcommand{\ed}{\end{defi}}
\newcommand{\er}{\end{rem}}
\newcommand{\ee}{\end{equation}}
\newcommand{\eea}{\end{eqnarray}}
\newcommand{\epr}{\end{proof}}
\newcommand{\ec}{\end{cor}}
\newcommand{\eco}{\end{cond}}
\newcommand{\eex}{\end{ex}}
\newcommand{\esch}{\end{sch}}

\newcommand{\nn}{ \nonumber \\ }
\newcommand{\pr}{{\em Proof:\ }}
\newcommand{\qed}{\vrule height 5pt width 5pt depth 0pt}
\newcommand{\ui}{^{(1)}}
\newcommand{\uii}{^{(2)}}
\newcommand{\di}{_{(1)}}
\newcommand{\dii}{{}_{(2)}}
\newcommand{\ot}{\otimes}

\newcommand{\ci}{\circ}
\newcommand{\lu}{\leftharpoonup}
\newcommand{\ld}{\leftharpoondown}
\newcommand{\ru}{\rightharpoonup}
\newcommand{\rd}{\rightharpoondown}

\newcommand{\End}{{\rm End}} 
\newcommand{\Hom}{{\rm Hom}} 
\newcommand{\M}{\mathcal{M}}
\newcommand{\calH}{\mathcal{H}}
\newcommand{\stac}[1]{\stackrel{\ot}{_{_{#1}}}}
\newcommand{\op}{^{op}}
\newcommand{\pri}{^{\prime}}
\newcommand{\inv}{^{-1}}
\newcommand{\h}{\ {\bf\underline{\ }}\ }
\newcommand{\il}{{\cal L}}

\newcommand{\lb}{\label}

\begin{document}

\large
\title{\bf Galois theory for Hopf algebroids}
 
\author{\sc Gabriella B\"ohm \\
Research Institute for Particle and Nuclear Physics, Budapest,\\
H-1525 Budapest 114, P.O.B. 49, Hungary\\
E-mail: G.Bohm@rmki.kfki.hu}
 
\date{}
 
\maketitle
\normalsize 
\begin{abstract} 
An extension $B\subset A$ of algebras over a commutative ring $k$ is an
$\calH$-extension for an $L$-bialgebroid $\calH$ if $A$ is an
$\calH$-comodule algebra and $B$ is the subalgebra of
its coinvariants. It is $\calH$-Galois if in addition the canonical map
$A\ot_B A\to A\ot_L \calH$ is an isomorphism or, equivalently, if the
canonical coring $(A\ot_L \calH:A)$ is a Galois coring.

In the case of a Hopf algebroid $\calH=(\calH_L,\calH_R,S)$, a comodule
algebra $A$ is defined as an algebra carrying compatible comodule structures
over both constituent bialgebroids $\calH_L$ and $\calH_R$.
If the antipode is bijective then $A$ is proven to be an $\calH_R$-Galois
extension of its coinvariants if and only if it is an $\calH_L$-Galois
extension. 

Results about bijective entwining structures
are extended to entwining structures over non-commutative
algebras in order to prove a Kreimer-Takeuchi type theorem 
for a finitely generated projective Hopf algebroid $\calH$ 
with a bijective antipode. It states that any $\calH$-Galois extension $B\subset
A$ is projective, and if $A$ is $k$-flat then already the
surjectivity of the canonical map implies the Galois property. 

The Morita theory, developed for corings by Caenepeel, Vercruysse and
Wang, is applied to obtain equivalent criteria for the Galois property of
Hopf algebroid extensions. This leads to Hopf algebroid
analogues of results for Hopf algebra extensions by Doi  and, in the
case of Frobenius Hopf algebroids, by Cohen, Fishman and Montgomery.
\end{abstract}

\section{Introduction}

An extension $B\subset A$ of algebras over a commutative ring $k$ is an
$H$-extension for a $k$-bialgebra $H$ if $A$ is a right
$H$-comodule algebra and $B$ is the subalgebra of
its coinvariants i.e. of elements $b\in A$ such that 
$b_{\langle 0\rangle}\ot b_{\langle 1\rangle} = b\ot 1_H$  
-- where the map $A\to A\ot_k H$, $a\mapsto  a_{\langle 0\rangle} \ot  
a_{\langle 1\rangle}$ is the coaction of 
$H$ on $A$ (summation understood). An $H$-extension $B\subset A$ is
$H$-Galois if the canonical map 
$A\ot_B A\to A\ot_k H$, $a\ot a\pri\mapsto aa\pri_{\langle 0\rangle} \ot 
a\pri_{\langle 1\rangle}$ is an isomorphism of $k$-modules.

In many cases it is technically much easier to check the surjectivity of the
canonical map than its injectivity. A powerful tool in the study of
$H$-extensions is the Kreimer-Takeuchi theorem \cite{KreiTak} stating that if
$H$ is a finitely generated projective Hopf algebra then the surjectivity of
the canonical map implies its bijectivity and also the fact that $A$ is
projective both as a left and as a right $B$-module.

The proof of the Kreimer-Takeuchi theorem went through both simplification and
generalization in the papers
\cite{Schau:Hopf-Gal,Schneider,Brze:Gal,SchauSchne}.  
In the present paper we adopt the method of Brzezi\'nski
\cite{Brze:Gal} and of Schauenburg and Schneider
\cite{SchauSchne}, who used the following observation. A comodule
algebra $A$ for a bialgebra $H$ determines a canonical entwining 
structure \cite{BrzeMaj} consisting of the algebra $A$, the coalgebra
underlying the bialgebra $H$, and the entwining map 
$H\ot_k A\to A\ot_k H$, $h\ot  a \mapsto a_{\langle 0\rangle} \ot h
a_{\langle 1\rangle}$. 
In the case when the bialgebra $H$ possesses a skew antipode,
this entwining map is a bijection. The proof of (a wide generalization of) the
Kreimer-Takeuchi theorem both in \cite{Brze:Gal} and in
\cite{SchauSchne} is based on the study of 
bijective entwining structures, under slightly different assumptions.
In Section \ref{sec:KT} below we show that these
arguments can be repeated almost without modification by using entwining
structures over 
non-commutative algebras \cite{internal}. 

In the paper \cite{Doi} Doi constructed a Morita context for an $H$-extension
$B\subset A$. 
If $H$ is finitely generated and projective as a $k$-module then the
surjectivity of one of the connecting maps is
equivalent to the projectivity and the Galois property of the
extension $B\subset A$, while the strictness of the Morita context is
equivalent to faithful flatness and the Galois property.
This observation made it possible to use all results of Morita theory for
characterizing $H$-Galois extensions. In the case when $H$ is a finite
dimensional Hopf algebra over a field (or a Frobenius Hopf algebra over a
commutative ring), the Morita context of Doi is equivalent to another Morita
context, introduced by Cohen, Fishman and Montgomery \cite{CohFishMont}.

One of the most beautiful applications of the theory of corings \cite{BrzeWis}
is the observation \cite{Brze} that the Galois property of an $H$-extension 
$B\subset A$ is
equivalent to the Galois property of a canonical $A$-coring $A\ot H$.
In \cite{CaenVerWang} the construction of the Morita context by  Doi has been
extended to any $A$-coring $C$ possessing a grouplike element (i.e. such
that $A$ is a $C$-comodule). 
In the case when $C$ is a finitely generated projective $A$-module (or
an $A$-progenerator, see \cite{Caenepeel}) the application of Morita theory
yields then 
several equivalent criteria for the Galois property of the coring $C$ and the
projectivity (or faithful flatness) of $A$ as a module for the subalgebra
of coinvariants in $A$.
In the case when the $A$-dual algebra of the coring $C$ is a
Frobenius extension of $A$, also the Morita context in \cite{CohFishMont} has
been generalized to the general setting of corings and the precise relation of
the two Morita contexts has been explained. 

The notion of bialgebra extensions has been generalized to bialgebroids by
Kadison \cite{Kadison} as an extension $B\subset A$ of $k$-algebras such that
$A$ is a comodule algebra and $B$ is the subalgebra of coinvariants. The
Galois property of a bialgebroid extension can be formulated also as the
Galois property of a canonical coring. This implies that the general theory,
developed in \cite{CaenVerWang}, can be applied also to bialgebroid extensions.

In the present paper we study Hopf algebroid extensions. The notion of Hopf
algebroids has been introduced in \cite{hgd, Brussels} and studied further in
\cite{integral}. It consists of two compatible (left and
right) bialgebroid structures on the same algebra which are related by the
antipode. 

In Section \ref{sec:Hopfext} 
a comodule of a Hopf algebroid is defined as a pair of comodules for both
constituent bialgebroids $\calH_L$ and $\calH_R$, 
in such a way that the $\calH_R$-coaction is $\calH_L$-colinear and the
$\calH_L$-coaction is $\calH_R$-colinear. 
In particular, a comodule algebra $A$ of a Hopf algebroid is defined as an
algebra carrying the structure of a 
compatible pair of comodule algebras of the constituent bialgebroids.
If the antipode of a Hopf algebroid $\calH$ is bijective, then we prove
that the $\calH_R$- and the $\calH_L$-coinvariants of any $\calH$-comodule
algebra $A$ coincide.  
What is more, we show that in this case $A$ is a Galois extension of its
coinvariant subalgebra by $\calH_L$
if and only if it is a Galois extension by $\calH_R$.

In Section \ref{sec:KT} it is shown that -- just as in the case of Hopf
algebras -- if $\calH$ is a Hopf algebroid with a bijective antipode then the
canonical entwining structure (over the non-commutative base algebra of
$\calH$), associated to an $\calH$-comodule algebra, is bijective. This fact
is used to prove a Kreimer-Takeuchi type theorem.

In Section \ref{sec:Morita} we apply the Morita theory for corings to a Hopf
algebroid extension $B\subset A$, looked at as an
extension by the constituent right bialgebroid $\calH_R$. 
In the finitely generated projective case this results in equivalent criteria,
under which $B\subset A$ is a projective $\calH_R$-Galois extension. 
Similarly, we can look at $B\subset A$ as an
extension by the constituent left bialgebroid $\calH_L$ and obtain equivalent
conditions for its projectivity and $\calH_L$ -Galois property.
Making use of the results about Hopf algebroid extensions in Section
\ref{sec:Hopfext}, and the Kreimer-Takeuchi type theorem proven in Section
\ref{sec:KT}, we conclude that if $\calH$ is a finitely generated projective
Hopf algebroid with a bijective antipode then the two equivalent sets of
conditions are equivalent also to each other.
In the case of Frobenius Hopf algebroids \cite{integral} we obtain a direct
generalization of (\cite{CohFishMont}, Theorem 1.2).

\smallskip

Throughout the paper $k$ is a commutative ring. By an algebra $R=(R,\mu,\eta)$
we 
mean an associative unital $k$-algebra. Instead of the unit map
$\eta$ we use sometimes the unit element $1_R\colon =\eta(1_k)$. We denote by
${_R\M}$,  
$\M_R$ and ${_R\M_R}$ the categories of left, right, and bimodules for
$R$, respectively. For the $k$-module of morphisms in ${_R\M}$,
${\M_R}$ and ${_R\M_R}$ we write ${_R\Hom(\ ,\ )}$, $\Hom_R(\ ,\ )$ and
${_R\Hom_R}(\ ,\ )$, respectively. 

\smallskip


{\bf Acknowledgment.} I would like to thank Tomasz Brzezi\'nski and
Korn\'el Szlach\'anyi for useful discussions.  
This work was supported by the Hungarian Scientific
Research Fund OTKA -- T 034 512 and T 043 159.

\section{Preliminaries}

\subsection{Bialgebroids and Hopf algebroids}
\lb{ss:hgd}

$L$ -bialgebroids \cite{Lu,Xu,Szlachanyi,Schauenburg} or, what were
shown in \cite{BrzeMi} to be equivalent to them, $\times_L$ -bialgebras 
\cite{Takeuchi} are generalizations of bialgebras to the case of
non-commutative base algebras. This means that instead of coalgebras and
algebras over commutative rings, one works with corings and rings 
over non-commutative base algebras. Recall that a coring over a
$k$-algebra $L$ is a comonoid in ${_L\M_L}$ while an $L$-ring is
a monoid in ${_L\M_L}$. The notion of $L$-rings is equivalent
to a pair, consisting of a $k$-algebra $A$ and an algebra homomorphism
$L\to A$.
\bd \label{lbgd} 
A {\em left bialgebroid} is a 6-tuple $\calH=(H,L,s,t,\gamma,\pi)$,
where $H$ and $L$ are $k$-algebras. $H$ is an $L\stac{k}
L\op$-ring via the algebra homomorphisms $s:L\to H$ and $t:L\op\to
H$, the images of which are required to commute in $H$. In terms of the maps
$s$ and $t$ one equips $H$ with an $L$-$L$ bimodule structure as
\be  l\cdot h\cdot l\pri\colon = s(l) t(l\pri)h \qquad \textrm{for}\ h\in H,\
l,l\pri\in L.\lb{Lbim}\ee
The triple $(H,\gamma,\pi)$ is an $L$-coring with respect to the bimodule
structure \eqref{Lbim}. Introducing Sweedler's
notation $\gamma(h)=h\di\stac{L} h\dii$ for 
$h\in H$ (with implicit summation understood), the axioms 
\bea  h\di t(l) \stac{L} h\dii &\!\!\!\!\!=&\!\!\!\!\! h\di \stac{L} h\dii
s(l) \lb{cros}\\ 
      \gamma(1_H)&\!\!\!\!\!=&\!\!\!\!\! 1_H\stac{L} 1_H \\
      \gamma(h h\pri)&\!\!\!\!\!=&\!\!\!\!\!\gamma(h) \gamma(h\pri) \lb{gmp} \\
      \pi(1_H) &\!\!\!\!\!=&\!\!\!\!\! 1_L \\
      \pi\left(h \ s\!\ci\!\pi(h\pri)\right)=&\!\!\!\!\!\pi(hh\pri)&\!\!\!\!=
      \pi\left(h \ t\!\ci\!\pi(h\pri)\right)\lb{aug}
\eea
are  required for all $l\in L$ and $h,h\pri\in H$.
\ed
Notice that -- although $H\stac{L} H$ is not an algebra -- axiom
(\ref{gmp}) makes sense in view of (\ref{cros}).

The bimodule (\ref{Lbim}) is defined in terms of multiplication by $s$
and $t$ on the left. The $R$-$R$ bimodule structure in a {\em right
bialgebroid} $\calH=(H,R,s,t,\gamma,\pi)$ is defined in terms of
multiplication on the right. For the details we refer to
\cite{KadSzlach}. 

The {\em opposite} of a left bialgebroid
$\calH=(H,L,s,t,\gamma,\pi)$ is the right bialgebroid
$\calH^{op}=(H^{op},L,t,s,$ $\gamma,\pi)$ where $H\op$ is the algebra
opposite to $H$. The {\em co-opposite} of $\calH$ is the left bialgebroid 
$\calH_{cop}=(H,L^{op},t,s,\gamma^{op},\pi)$ where $\gamma^{op}:H\to
H\stac{L^{op}} H$ is the opposite coproduct $h\mapsto
h\dii\stac{L^{op}} h\di$.

It has been observed in \cite{KadSzlach} that for a left bialgebroid
$\calH=(H,L,s,t,\gamma,\pi)$, such that $H$ is finitely generated and
projective as a left or right  $L$-module, the corresponding $L$-{\em
dual} possesses a right bialgebroid structure over the base algebra
$L$. Analogously, the $R$-duals of a finitely generated projective 
right bialgebroid over $R$ possess left
bialgebroid structures. For the explicit forms of the dual bialgebroid
structures consult \cite{KadSzlach}.

Before defining the notion of a Hopf algebroid, let us introduce some
notations. Analogous notations were used already in
\cite{hgd,integral}.

When dealing with an $L\stac{k}L\op$-ring $H$, 
we have to face the situation that $H$ 
carries different module structures over the base algebra $L$. In this
situation the usual notation $H\stac{L} H$ would be ambiguous. Therefore we
make the following notational convention. In terms of the algebra homomorphisms
$s:L\to H$ and $t:L\op\to H$ (with commuting images in $H$) we introduce four
$L$-modules 
\bea {_L H}:&\qquad &l\cdot h \colon = s(l)h\nn
H_L:&\qquad &h\cdot l\colon = t(l)h\nn
H^L:&\qquad&h\cdot l=hs(l)\nn
{^L H}:&\qquad& l\cdot h= ht(l).
\lb{amod}\eea
This convention can be memorized as left indices stand for left modules and
right indices stand for right modules. Upper indices refer to modules defined
in terms of right multiplication and lower indices refer to the ones defined
in terms of left multiplication.

In writing $L$-module tensor products, we write out explicitly the module
structures of 
the factors that are taking part in the tensor products, and do not put marks
under the symbol $\ot$. E.g. we write $H_L\ot {_L H}$. In writing elements of
tensor product modules we do not distinguish between the various module
tensor products. That is, we write both $h\stac{L} h\pri\in H_L\ot {_L H}$ and
$g\stac{L} 
g\pri\in H^L\ot {_L H}$, for example.

A left $L$-module can be considered canonically as a right $L\op$-module, and
sometimes we want to take a module tensor product over $L\op$. In this case we
use the name of the corresponding $L$-module and the fact that the tensor
product is taken over $L\op$ should be clear from the order of the factors. 

In writing multiple tensor products we use different types of letters to
denote which module structures take part in the same tensor product. 
\bd \lb{hgd}
A {\em Hopf algebroid} $\calH=(\calH_L,\calH_R,S)$ consists of a 
left bialgebroid $\calH_L=(H,L,s_L,$ $t_L,\gamma_L,\pi_L)$ and a right 
bialgebroid $\calH_R=(H,R,s_R,t_R,\gamma_R,\pi_R)$, with common total algebra
$H$, and a $k$-module 
map $S:H\to H$, called the antipode, such that the following axioms
hold true
\footnote[1]{An equivalent set of axioms is obtained by requiring $S$
    to be only an $R$-$L$ bimodule map, see Remark 2.1 in \cite{BB:cleft}.}: 
\bea i)&& s_L\ci \pi_L\ci t_R=t_R,\qquad t_L\ci \pi_L\ci s_R=s_R
\quad \textrm{and}\nn
&& s_R\ci \pi_R\ci t_L=t_L,\qquad t_R\ci \pi_R\ci s_L=s_L;
\lb{hgdi}\\
ii)&&(\gamma_L\ot {^R H})\ci \gamma_R=(H_L\ot \gamma_R)\ci\gamma_L
\quad \textrm{as\ maps\ } H\to H_L\ot {_L H^R}\ot {^R H}\quad
\textrm{and}\nn
&&(\gamma_R\ot {_L H})\ci \gamma_L=(H^R\ot \gamma_L)\ci\gamma_R
\quad \textrm{as\ maps\ } H\to H^R\ot {^R H_L}\ot {_L H};\lb{hgdii}\\
iii)&& S\ \textrm{is\ both\ an\ }L\textrm{-}L\ \textrm{bimodule\ map\ } {^L
H_L}\to {_L H^L}\ \textrm{and\ an\ } R\textrm{-}R\ \textrm{bimodule\ map}\nn
&& {^R H_R}\to {_R H^R}; \lb{hgdiii}\\
iv)&&\mu_H\ci (S\ot {_L H})\ci \gamma_L=s_R\ci\pi_R\quad \textrm{and}\nn
   &&\mu_H\ci (H^R\ot S)\ci \gamma_R=s_L\ci \pi_L. \lb{hgdiv}
\eea
\ed
If $\calH=(\calH_L,\calH_R,S)$ is a Hopf algebroid then so is
$\calH\op_{cop}=((\calH_R)\op_{cop},(\calH_L)\op_{cop},S)$ and if $S$ is
bijective then also $\calH_{cop}=((\calH_L)_{cop},(\calH_R)_{cop},S\inv)$ and
$\calH\op=((\calH_R)\op,(\calH_L)\op,S\inv)$. 

We are going to use the following variant of the Sweedler-Heynemann index
convention. For a Hopf algebroid $\calH=(\calH_L,\calH_R,S)$ we 
use the notation $\gamma_L(h)=h\di\ot h\dii$ with lower indices and
$\gamma_R(h)=h\ui\ot h\uii$ with upper indices for $h\in H$ in the
case of the coproducts of $\calH_L$ and 
$\calH_R$, respectively. In both cases implicit summation is understood.
Axioms (\ref{hgdii}) read in this notation as 
\bea {h\ui}\di\ot {h\ui}\dii\ot h\uii&=&h\di\ot {h\dii}\ui\ot
{h\dii}\uii\nn
{h\di}\ui\ot {h\di}\uii\ot h\dii&=&h\ui\ot {h\uii}\di \ot {h\uii}\dii
\nonumber \eea
for $h\in H$.

It is proven in (\cite{integral}, Proposition 2.3) that the base algebras 
$L$ and $R$ of the left and right bialgebroids $\calH_L$ and $\calH_R$ 
in a Hopf algebroid $\calH$ are anti-isomorphic via
any of the isomorphisms $\pi_L\ci s_R$ and $\pi_L\ci t_R$. The antipode is a
homomorphism of left bialgebroids $\calH_L\to (\calH_R)^{op}_{cop}$ and also
$(\calH_R)^{op}_{cop}\to \calH_L$ in the sense that it is an anti-algebra
endomorphism of $H$ and the pair of maps $(S,\pi_L\ci s_R)$ 
is a coring homomorphism from the $R\op$-coring $(H,\gamma_R^{op},\pi_R)$ to
the $L$-coring $ (H,\gamma_L,\pi_L)$
and the pair of algebra homomorphisms $(S,\pi_R\ci s_L)$ is a coring
homomorphism $(H,\gamma_L,\pi_L)\to (H,\gamma_R^{op},\pi_R)$.

We term a Hopf algebroid $\calH$, for that all
modules $H^R$, ${^R H}$, $H_L$ and ${_L H}$ is finitely generated and
projective, as a {\em finitely generated projective} Hopf algebroid.

Left {\em integrals} in a left bialgebroid $\calH_L$ are defined
(\cite{hgd}, Definition 5.1) as the invariants of the left regular $H$-module
i.e. the elements of 
$$\il(H)\colon=\{\ \ell\in H\ \vert\ h\ell=s_L\!\ci\!\pi_L(h)\ \ell\quad
\forall h\in H\ \}.$$
By (\cite{integral}, Scholium 2.8) an element $\ell$ of a Hopf algebroid
$\calH$ is a left integral if and only if $h\ell\ui\stac{R} S(\ell\uii)=
\ell\ui\stac{R} S(\ell\uii)h$ for all $h\in H$.

A left integral $\ell$ in a Hopf algebroid $\calH$ is called {\em
non-degenerate} (\cite{hgd}, Definition 5.3) if both maps
\bea 
\ell_R:H^*\colon = \Hom_R(H^R,R) \to H&\qquad &\phi^*\mapsto \phi^*\ru
\ell\equiv \ell\uii \ t_R\!\ci\! \phi^*(\ell\ui)\quad \textrm{and}\nn
{_R\ell}:{^*H}\colon = {_R\Hom}({^RH},R) \to H&\qquad &{^*\phi}\mapsto
{^*\phi}\rd \ell\equiv \ell\ui \ s_R\!\ci\! {^*\phi}(\ell\uii)
\nonumber\eea
are isomorphisms. By (\cite{hgd}, Proposition 5.10) for a non-degenerate left
integral $\ell$ in a Hopf algebroid $\calH$ also the maps 
\bea 
\ell_L:H_*\colon = \Hom_L(H_L,L) \to H&\qquad &\phi_*\mapsto \ell\lu \phi_*
\equiv s_L\!\ci\! \phi_*(\ell\di)\ \ell\dii \quad \textrm{and}\nn
{_L\ell}:{_*H}\colon = {_L\Hom}({_LH},L) \to H&\qquad &{_*\phi}\mapsto
\ell\ld {_*\phi}\equiv t_L\!\ci\! {_*\phi}(\ell\dii) \ \ell\di
\nonumber\eea
are isomorphisms.
It is shown in (\cite{integral}, Theorem 4.7, 
see also the {\em Corrigendum}) that the
existence of a non-degenerate left integral in a 
finitely generated projective
Hopf algebroid $\calH$ is
equivalent 
to the Frobenius property of any of the four extensions $s_R:R\to H$,
$t_R:R\op\to H$, $s_L:L\to H$ and $t_L:L\op\to H$ and it implies the
bijectivity of the antipode. What is more, if the Hopf algebroid $\calH$
possesses a non-degenerate left integral then also the four duals $H^*$,
${^*H}$, $H_*$ and ${_* H}$ possess (anti-) isomorphic Hopf algebroid
structures with non-degenerate integrals (\cite{hgd}, Theorem 5.17 and
Proposition 5.19). Motivated by these results we term a Hopf algebroid
possessing a non-degenerate left integral as a {\em Frobenius} Hopf algebroid.
Recall from \cite{Szla:Doubalg} that a Frobenius Hopf
algebroid is equivalent to a distributive Frobenius double algebra.

\subsection{Module and comodule algebras}
\lb{ss:comalg}

The category ${_H\M}$ of left modules for the total algebra $H$ of a left 
bialgebroid $(H,L,s,t,\gamma,\pi)$ is a monoidal category. As a matter of fact,
any $H$-module is an $L$-$L$ bimodule via $s$ and $t$. The monoidal product in
${_H \M}$ is the $L$-module tensor product with $H$-module structure 
$$h\cdot (m\stac{L} n)\colon =h\di\cdot m \stac{L} h\dii\cdot n\qquad
\textrm{for}\ h\in H, \ m\stac{L} n\in M\stac{L} N$$
and the monoidal unit is $L$ with $H$-module structure 
$$ h\cdot l\colon =\pi(hs(l))\qquad \textrm{for}\ h\in H, l\in L.$$
A left {\em $H$-module algebra} is defined as a monoid in the monoidal
category ${_H \M}$.
A left $H$ -module algebra $A$ is in particular an
$L$-ring via the homomorphism
$$L\to A\qquad l\mapsto l\cdot 1_A\equiv 1_A\cdot l.$$
The {\em invariants} of $A$ are the elements of
$$A^H\colon =\{\ a\in A\ \vert\ h\cdot a=s\!\ci\!\pi(h)\cdot a\quad \forall
h\in H\ \}.$$
Just in the same way, the category $\M_H$ of right modules for the total algebra
$H$ of a right $R$-bialgebroid is a monoidal category with monoidal
product the $R$-module tensor product and monoidal unit $R$. A right
$H$-module algebra is a monoid in $\M_H$.
A right module  algebra is in particular
an $R$-ring. The invariants are defined analogously to the left case in terms
of the counit. 

\bigskip

By a comodule for a left bialgebroid $\calH=(H,L,s,t,\gamma,\pi)$ we mean a
comodule for the $L$-coring $(H,\gamma,\pi)$. Recall that the 
category of left $\calH$-comodules is also a monoidal category in the
following way. Any left $\calH$-comodule $(M,\tau)$ can be equipped with a
right $L$-module structure via  
\be m\cdot l\colon = \pi(m_{\langle -1\rangle} s(l))\cdot m_{\langle 0 \rangle}
\qquad \textrm{for}\ m\in M,\ l\in L\lb{rightac}\ee
where $m_{\langle -1\rangle}\stac{L} m_{\langle 0\rangle}$ stands for
$\tau(m)$ (summation understood). Indeed, (\ref{rightac}) is the
unique right action via which $M$ becomes an $L$-$L$
bimodule and $\tau$ becomes an $L$-$L$ bimodule map from ${_L M_L}$ to
the Takeuchi product $H\times_L M$. Recall from \cite{Takeuchi} that $H\times_L
M$ is the $L$-$L$ submodule of ${_L H^L_{\bf L}}\ot {_{\bf L} M}$ the elements
$\sum_i h_i\stac{L} m_i$ of which satisfy  
\be
\sum_i h_i\stac{L} m_i\cdot l=  \sum_i h_i t(l)\stac{L} m_i
\qquad \textrm{for}\ \ l\in L.
\ee
This observation amounts to saying that our definition of left 
$\calH$-comodules is equivalent to (\cite{Schauenburg}, Definition
5.5). Hence, without loss of generality, from now on we can
think of a 
comodule in this latter sense. On the basis of (\cite{Schauenburg}, Definition
5.5) the category $^H\M$ of left $\calH$ -comodules was shown in
(\cite{Schauenburg}, Proposition 5.6) to 
be monoidal. The monoidal product is the $L$-module tensor product
with comodule structure 
$$M\stac{L} N\to H\stac{L} M\stac{L} N\qquad m\stac{L} n\mapsto 
m_{\langle -1\rangle} n_{\langle -1\rangle}\stac{L} m_{\langle 0\rangle} 
\stac{L} n_{\langle 0\rangle}$$
and the monoidal unit is $L$ with comodule structure 
$$L\to L\stac{L} H\simeq H\qquad l\mapsto s(l).$$
Following (\cite{Schauenburg}, Definition 5.7), a left {\em comodule
algebra} for a left bialgebroid $\calH$ is a monoid in the monoidal category
${^H\M}$. Notice that -- in view of the equivalence of the two definitions
of $\calH$ -comodules --  this definition of comodule algebras is
equivalent to (\cite{BrzeCaenMil}, Definition 3.4).

By similar arguments also the category of right $\calH$-comodules -- that is
of right comodules for the $L$-coring $(H,\gamma,\pi)$ -- is
monoidal. The 
monoidal product is the $L$-module tensor product with coaction
$$M\stac{L} N\to M\stac{L} N \stac{L} H\qquad m\stac{L} n\mapsto 
m_{\langle 0\rangle}\stac{L}n_{\langle 0\rangle}\stac{L}n_{\langle 1\rangle}
m_{\langle 1\rangle}$$
and the monoidal unit is $L$ with comodule structure
$$L\to H\qquad l\mapsto t(l),$$
where the left $L$-module structure of a right $(H,\gamma,\pi)$-comodule $M$ is
defined as $l\cdot m\colon = m_{\langle 0\rangle}\cdot \pi(m_{\langle
  1\rangle} s(l))$.
Similarly to the case of left comodule algebras we expect the coaction of a
right comodule algebra $A$ to be multiplicative
-- i.e. such that $(aa\pri)_{\langle 0\rangle}\stac{L} 
(aa\pri)_{\langle 1\rangle}=a_{\langle 0\rangle}a\pri_{\langle 0\rangle}
\stac{L}a_{\langle 1\rangle}a\pri_{\langle 1\rangle}$ for $a,a\pri\in A$.
Therefore we consider the monoidal category $(\M^H)\op$, the monoidal
structure of which is the opposite of $\M^H$ (i.e. it comes from the monoidal
structure of ${_{L\op}\M_{L\op}}$). A right $\calH$-comodule algebra is defined
as a monoid in the category $(\M^H)\op$.

Notice that a left $\calH$-comodule algebra is in particular an
$L$-ring while a right $\calH$-comodule algebra is an $L\op$-ring.

The {\em coinvariants} of a left (right) $\calH$-comodule algebra $(A,\tau)$
are the elements of the subalgebra
$$ A^{coH}\colon = \{\ a\in A\ \vert\ \tau(a)=1_H\stac{L} a\ \}\quad
(\ A^{coH}\colon = \{\ a\in A\ \vert\ \tau(a)=a\stac{L} 1_H\ \}\ ).
$$

Recall that for a left $L$-bialgebroid $\calH$ the left and right
$L$-duals ${_*H}$ and $H_*$ are rings. 
There is a faithful functor from the category  $\M^H$ of right
$\calH$-comodules to the category $\M_{{}_*H}$ of right ${_*
H}$-modules which is an isomorphism if and only if the
module $_LH$ is finitely generated and projective.
There is a faithful functor also from ${^H \M}$ to $\M_{H_*}$ 
 which is an isomorphism if and only if the module $H_L$ is finitely generated
 and projective. 

Left and right comodules, comodule algebras and their coinvariants for a
right bialgebroid are defined analogously. 

\subsection{Entwining structures over non-commutative algebras}
\lb{ss:entw}

Entwining structures over non-commutative algebras were introduced in
\cite{internal} as mixed distributive laws in the bicategory of [Algebras,
Bimodules, Bimodule maps]. This definition is clearly equivalent to a monad in
the bicategory of corings i.e. the bicategory of comonads \cite{Street} in the 
bicategory of [Algebras, Bimodules, Bimodule maps]. Explicitly, we have
\bd An {\em entwining structure} over an algebra $R$ is a triple $(A,C,\psi)$
where $A=(A,\mu,\eta)$ is an $R$-ring, $C=(C,\Delta,\epsilon)$ is an
$R$-coring and $\psi$ is an $R$-$R$ bimodule map $C\stac{R} A\to A\stac{R} C$
satisfying 
\bea 
&& \psi\ci(C\stac{R} \eta)=\eta\stac{R} C\nn
&&(A\stac{R} \epsilon)\ci \psi=\epsilon\stac{R} A\nn
&&(\mu\stac{R} C)\ci (A\stac{R} \psi)\ci(\psi\stac{R}
A)=\psi\ci(C\stac{R}\mu)\nn
&&(A\stac{R} \Delta)\ci \psi=(\psi\stac{R} C)\ci (C\stac{R} \psi)\ci
(\Delta\stac{R} A).
\nonumber\eea
An entwining structure is {\em bijective} if $\psi$ is an isomorphism.
\ed
It is shown in (\cite{internal}, Example 4.5) that an entwining structure
$(A,C,\psi)$ over the algebra $R$ determines an $A$-coring structure on
$A\stac{R} C$ with $A$-$A$ bimodule structure
$$a_1\cdot (a\stac{R} c)\cdot a_2= a_1 a \psi(c\stac{R} a_2)\qquad
\textrm{for}\ a_1,a_2\in A,\  a\stac{R} c\in A\stac{R} C,$$
coproduct $A\stac{R} \Delta$ and counit $A\stac{R}\epsilon$.
\bd A right {\em entwined module} over an $R$-entwining structure $(A,C,\psi)$
is a right comodule over the corresponding $A$-coring $A\stac{R}
C$. Explicitly, it 
is a triple $(M,\rho,\tau)$, where $(M,\rho)$ is a right $A$-module, making $M$
in particular a right $R$-module. The pair $(M,\tau)$ is a right $C$-comodule
such that $\tau$ is a right $A$-module map i.e. 
$$\tau\ci \rho=(\rho\stac{R} C)\ci (M\stac{R} \psi)\ci (\tau\stac{R} A). $$
A {\em morphism} of entwined modules is a morphism of comodules for the
$A$-coring $A\stac{R} C$, that is, an $A$-linear and $C$-colinear map.
The category of entwined modules will be denoted by $\M_A^C$.

The {\em coinvariants} of an entwined module are its coinvariants for the
$A$-coring $A\stac{R} C$. If the $R$-coring $C$ possesses a grouplike element,
then this the same as $C$-coinvariants.
\ed

By (\cite{BrzeWis}, 18.13 (2)) the forgetful functor $\M^C_A\to \M_A$ possesses
a right adjoint, the functor 
\be \h \stac{R} C:\M_A\to \M_A^C\qquad (\ M,\rho\ )\mapsto (\ M\stac{R} C\ ,\  
(\rho\stac{R}C)\ci (M\stac{R}\psi)\ ,\  M\stac{R} \Delta\ ),\lb{adjC}\ee
where the right $R$-module structure of $M$ comes from its $A$-module
structure.
What is more, by the self-duality of the notion of $R$-entwining
structures, also (\cite{BrzeWis}, 32.8 (3)) extends to entwining structures
over non-commutative algebras. That is, also the forgetful functor $\M^C_A\to
\M^C$ possesses a left adjoint, the functor
\be \h\stac{R} A:\M^C\to \M^C_A\qquad (\ M,\tau\ )\mapsto (\ M\stac{R} A\ ,\ 
(M\stac{R}\mu)\ ,\  (M\stac{R}\psi)\ci(\tau\stac{R} A)\ ).\lb{adjA}\ee
This implies, in particular, that both $A\stac{R} C$ and $C\stac{R} A$ are
entwined modules. The morphism $\psi$ becomes a morphism of entwined modules. 

If the $R$-coring $C$ possesses a grouplike element $e$, then also the
$A$-coring $A\stac{R} C$ possesses a grouplike element $1_A\stac{R} e$. Hence
$A$ is an entwined module via the right regular $A$-action and the
$C$-coaction 
$$A\to A\stac{R} C\qquad a\mapsto \psi(e\stac{R} a).$$
In this case also the functors 
\be (\h)^{coC}:\M_A^C\to \M_{A^{coC}}\qquad \textrm{and}\qquad
\h\stac{A^{coC}} A:\M_{A^{coC}}\to \M_A^C\lb{adj}\ee
are adjoints (\cite{BrzeWis}, 28.8). (The entwined module structure of
$N\stac{A^{coC}} A$, for a right $A^{coC}$-module $N$, is defined via the second
tensor factor). The unit and the counit of the adjunction are
\bea
&\eta_N: N\to (N\stac{A^{coC}} A)^{coC}\qquad & n\mapsto n\stac{A^{coC}}
1_A\quad \textrm{and}\nn
& \mu_M: M^{coC}\stac{A^{coC}} A\to M\qquad &m\stac{A^{coC}} a\mapsto m\cdot a
\nonumber\eea
for any right $A^{coC}$-module $N$ and entwined module $M$.

\subsection{Morita theory for corings}
\lb{ss:Morcor}

In the paper \cite{CaenVerWang} a Morita context
$(A^{^*C},{^*C},A,{^*C}^{^*C},\nu,\mu)$ has been associated to an $A$-coring
$C$ possessing a grouplike element $e$. Here the ring
${^*C}={_A\Hom(C,A)}$ is the left $A$-dual 
of the $A$-coring $C$ with multiplication 
$(fg)(c)=g\left( c\di\cdot f(c\dii)\right)$. The invariants of a right
${^*C}$-module $M$ are defined with the help of the grouplike element $e$ as
the elements of 
$$M^{^*C}\colon =\{\ m\in M\ \vert\ m\cdot f=m\cdot [\epsilon(\h)f(e)]\quad
\forall f\in {^*C}\ \}.$$
In terms of the grouplike element $e$, the $k$-module $A$ can be equipped with
a right ${^*C}$-module structure as
$$ a\cdot f\colon = f(e\cdot a)\qquad \textrm{for}\ a\in A, f\in {^*C}.$$
The ring $A^{^*C}$ is the subring of ${^*C}$-invariants of $A$ i.e.
$$A^{^*C}=\{\ b\in A\ \vert\ f(e\cdot b)=bf(e)\quad \forall f\in {^*C}\ \}. $$
$A$ is an $A^{^*C}$ - ${^*C}$ bimodule via
$$ b\cdot a\cdot f\colon = bf(e\cdot a)=f(e\cdot ba)\qquad \textrm{for}\ b\in
  A^{^*C}, a\in A, f\in {^*C}.$$
${^*C}^{^*C}$ is the $k$-module of ${^*C}$-invariants of the right regular
${^*C}$-module i.e.
$$ {^*C}^{^*C}=\{\ q\in {^*C}\ \vert \ f(c\di\cdot q(c\dii))=f(q(c)\cdot
  e)\quad \forall f\in {^*C}, \ c\in C\ \}.$$
It is a ${^*C}$ - $A^{^*C}$ bimodule via 
$$(f\cdot q\cdot b)(c)\colon = q\left(c\di\cdot f(c\dii)\right)b\qquad
\textrm{for} f\in {^*C},\ q\in {^*C}^{^*C},\ b\in A^{^*C}, c\in C.$$
The connecting maps $\nu$ and $\mu$ are given as
\bea
&\nu:A\stac{^*C} {^*C}^{^*C}\to A^{^*C}\qquad &a\stac{^*C} q\mapsto q(e\cdot
a)\quad 
\textrm{and}\nn
&\mu:{^*C}^{^*C}\stac{A^{^*C}}A\to {^*C}\qquad &q\stac{A^{^*C}}a\mapsto
(\ c\mapsto q(c)a\ ).
\nonumber\eea
In \cite{CaenVerWang} the following theorem has been proven.
\bt \lb{CVW}
Let $C$ be an $A$-coring possessing a grouplike element $e$, and let  
$(A^{^*C},{^*C},A,{^*C}^{^*C},\nu,\mu)$ be the Morita context associated to
it.

{\em (1)} (\cite{CaenVerWang}, Theorem 3.5). If $C$ is finitely generated and
    projective as a left $A$-module (hence the categories $\M^C$ and
    $\M_{^*C}$ are isomorphic and the $C$-coinvariants coincide with the
    ${^*C}$-invariants) then the following assertions are equivalent.

{\em (a)} The map $\mu$ is surjective (and, a fortiori, bijective).

{\em (b)} The functor $(\h)^{coC}:\M^C\to \M_{A^{coC}}$ is fully faithful.

{\em (c)} $A$ is a right ${^*C}$-generator.

{\em (d)} $A$ is projective as a left $A^{coC}$-module and the map
$${^*C}\to {_{A^{coC}} \End} (A)\qquad f\mapsto \left(\ a\mapsto f(e\cdot a)\
      \right)$$
is an algebra anti-isomorphism.

{\em (e)}  $A$ is projective as a left $A^{coC}$-module and the $A$-coring $C$
with grouplike element $e$ is a Galois coring.

{\em (2)} (\cite{CaenVerWang}, Theorem 2.7). If the algebra extension
$$A\to {^*C}\qquad a\mapsto \left(\ c\mapsto \epsilon(c)a\ \right)$$
is a Frobenius extension with Frobenius system $(\psi, u_i\stac{A} v_i)$ then
the Morita context $(A^{^*C},{^*C},A,$ ${^*C}^{^*C},\nu,\mu)$ is equivalent to
the Morita context $(A^{^*C},{^*C},A,A,\nu\pri,\mu\pri)$ via the isomorphism
$$A\to {^*C}^{^*C}\qquad a\mapsto (\ c\mapsto \sum_i v_i[c\cdot au_i(e)]\ ).$$
\et


\section{Hopf algebroid extensions}
\lb{sec:Hopfext}

An algebra extension $B\subset A$ is an $\calH_R$ -extension for a right
bialgebroid $\calH_R=(H,R,s,t,\gamma,\pi)$ if $A$ is a right
$\calH_R$-comodule algebra and $B$ is the subalgebra of
$\calH_R$-coinvariants of $A$. In this situation -- denoting the
$R$-coring $(H,\gamma,\pi)$ by $C$ -- the triple $(A,\calH_R,C)$ is a
Doi-Koppinen datum over $R$ in the sense of (\cite{BrzeCaenMil},
Definition 3.6). This implies that the $R$-$R$ bimodule map 
\be\lb{entw}
\psi:H\stac{R} A\to A\stac{R} H\qquad h\stac{R} a \mapsto 
a^{\langle 0\rangle}\stac{R} ha^{\langle 1\rangle}\ee
gives rise to an entwining structure $(A,C,\psi)$ over $R$. Hence the
$R$-$R$ bimodule $A\stac{R} H$ possesses an $A$-coring structure with
left and right $A$-actions
$$a_1\cdot(a\stac{R} h)\cdot a_2 = 
a_1 a a_2^{\langle 0\rangle}\stac{R} ha_2^{\langle 1\rangle}
\qquad \textrm{for}\ a_1,a_2\in A,\ a\stac{R} h\in A\stac{R} H,$$
coproduct $A\stac{R} \gamma$ and counit $A\stac{R}\pi$. This coring
possesses a grouplike element $1_A\stac{R} 1_H$. The $\calH_R$-extension
$B\subset A$ was termed $\calH_R$-Galois in \cite{Kadison} if the
$A$-coring $A\stac{R} H$, associated to it above, is a Galois
coring. This 
means bijectivity of the canonical map
\be{\rm can}_R: A\stac{B} A\to A\stac{R} H\qquad a\stac{B} a\pri \mapsto
a a^{\prime\langle 0\rangle} \stac{R}  a^{\prime\langle 1\rangle}.\lb{canR}\ee
Analogously, in the case of a right comodule algebra $A$ for the left
bialgebroid $\calH_L=(H,L,s,t,\gamma,$ $\pi)$ the $\calH_L$-Galois
property of the extension 
$A^{co\calH_L}\subset A$ means the bijectivity of the canonical map
\be {\rm can}_L: A\stac{A^{co\calH_L}} A\to A\stac{L} H\qquad
a\stac{A^{co\calH_L}} a\pri 
\mapsto 
a_{\langle 0\rangle} a\pri \stac{L}  a_{\langle 1\rangle}.\lb{canL}\ee
This is equivalent to the Galois property of the $A$-coring $A\stac{L}
H$ with $A$-$A$ bimodule structure
$$a_1\cdot(a\stac{L} h)\cdot a_2 = 
a_{1\langle 0\rangle} a a_2\stac{L} a_{1\langle 1\rangle}h
\qquad \textrm{for}\ a_1,a_2\in A,\ a\stac{L} h\in A\stac{L} H,$$
coproduct $A\stac{L} \gamma$, counit $A\stac{L}\pi$  and grouplike
element $1_A\stac{L} 1_H$. (Recall from Section \ref{ss:comalg}
that by our convention $A$ is an $L\op$-ring in this case.)

\bp\lb{prop:withdrawn}
{\color{red}
This is withdrawn because of an unjustified step in the proof.}
\ep

Withdrawal of Proposition \ref{prop:withdrawn} means that -- although we are
not aware of any counterexample -- it is {\em no longer proven} that for a
Hopf algebroid $\calH=(\calH_L,\calH_R,S)$, any $\calH_L$-comodule possesses
an $\calH_R$-comodule structure, or vice versa. 
(It is discussed in \cite[(arXiv version) Theorem 2.5]{BB:cleft} under what
additional assumptions this can be proven.)
Therefore, 
{\color{red}
Definition \ref{hopfcomod} and Lemma \ref{LReq} need to be modified}
as follows. No other results in the rest of the paper are affected by these
changes. 

\bd \lb{hopfcomod}
A right {\em comodule for the Hopf algebroid} $\calH=(\calH_L,\calH_R,S)$ is
a triple $(M,\tau_L,\tau_R)$, where the pair $(M,\tau_R)$ is a right
$\calH_R$-comodule and $(M,\tau_L)$ is a right $\calH_L$-comodule such that
the $R$-$R$ and the $L$-$L$ bimodule structures of $M$ are related via
\be
l\cdot m\cdot l'=\pi_R\ci t_L(l')\cdot m\cdot \pi_R\ci t_L(l),\label{eq:lbim}
\ee
and the compatibility relations 
\bea
&& (\tau_R\stac{L} H)\ci \tau_L=(M\stac{R} \gamma_L)\ci \tau_R\lb{eq:trcolin}\\
&& (\tau_L\stac{R} H)\ci \tau_R=(M\stac{L} \gamma_R)\ci \tau_L\lb{eq:tlcolin}
\eea 
hold true. 

The right $\calH$-comodule $(A,\tau_L,\tau_R)$ is said to be a right {\em
  $\calH$-comodule algebra} if $(A,\tau_R)$ is a right
$\calH_R$-comodule algebra and $(A,\tau_L)$ is a right
$\calH_L$-comodule algebra. 
\footnote[2]{It is proven in the arXiv version and the Corrigendum of
  \cite{BB:cleft} that the category ${\mathcal M}^\calH$ of right comodules of
  a Hopf algebroid $\calH$ is a monoidal category, with strict monoidal
  forgetful functors to the comodule categories of the constituent
  bialgebroids. Consequently, an $\calH$-comodule algebra is the same as a
  monoid in ${\mathcal M}^\calH$.} 
\ed
We follow the convention of using upper indices to denote the
components of the coaction of a right bialgebroid and lower indices
in the case of a left bialgebroid.
\bl \lb{LReq}
Let $\calH=(\calH_L,\calH_R,S)$ be a Hopf algebroid with a bijective antipode
and let $A$ be a right $\calH$-comodule algebra. 
Then the subalgebras of $\calH_R$-coinvariants and of $\calH_L$-coinvariants
in $A$ coincide. Moreover, denoting this coinvariant subalgebra by $B$,
the canonical map
$${\rm can}_R: A\stac{B} A\to A\stac{R} {H}\qquad a\stac{B} a\pri \mapsto
a a^{\prime\langle 0\rangle} \stac{R}  a^{\prime\langle 1\rangle}$$
is injective/surjective/bijective if and only if the canonical map
$${\rm can}_L: A\stac{B} A\to A\stac{L}{H}\qquad a\stac{B} a\pri \mapsto
a_{\langle 0\rangle} a\pri \stac{L}  a_{\langle 1\rangle}$$
is injective/surjective/bijective.
\el

For a Hopf algebroid  $\calH=(\calH_L,\calH_R,S)$ with a bijective antipode,
and an $\calH$-comodule 
algebra $A$ with $\calH_R$-, equivalently, $\calH_L$-coinvariant subalgebra
$B$, we term the algebra extension $B\subset A$ an {\em $\calH$-extension}.
  
In the case when both canonical maps in Lemma \ref{LReq} are bijective, we say
that $B\subset A$ is an {\em $\calH$-Galois extension}.
\smallskip

{\em Proof of Lemma \ref{LReq}:}
For any right $\calH$-comodule $(M,\tau_L,\tau_R)$, there exists an
isomorphism
$$\Phi_M:M\stac{R} {H}\to M\stac{L} {H}\qquad m\stac{R} h\mapsto 
m_{\langle 0\rangle}\stac{L} m_{\langle 1\rangle} S(h)$$
with inverse 
$$\Phi_M\inv:M\stac{L} {H}\to M\stac{R} {H}\qquad m\stac{L} h\mapsto 
m^{\langle 0\rangle}\stac{R} S\inv(h) m^{\langle 1\rangle}.$$
Since $\Phi_M(\tau_R(m))=m\ot_L 1_H$ and $\Phi_M(m\ot_R 1_H)=\tau_L(m)$,
it follows that in particular the $\calH_R$-coinvariants and the
$\calH_L$-coinvariants of any $\calH$-comodule algebra $A$ coincide.

Using the $\calH_R$-colinearity of $\tau_L$, the Hopf algebroid axiom
(\ref{hgdiv}), the fact that the image of $\tau_L$ is in the Takeuchi product
$A\times_L H$, the $L^{op}$-ring structure of $A$ and the counitality of
$\tau_L$, one checks that for an $\calH$-comodule algebra
$(A,\tau_L,\tau_R)$ the identity $\Phi_A\ci {\rm can}_R={\rm can}_L$ holds
true. 
\hfill\qed
\smallskip

\bex A $k$-bialgebra $(H,\mu,\eta,\Delta,\epsilon)$ is an example both of
a left- and of a right bialgebroid via the correspondence ${\cal B}\colon
=(H,k,\eta,\eta,\Delta,\epsilon)$. A Hopf algebra
$(H,\mu,\eta,\Delta,\epsilon,S)$ is an example of a Hopf algebroid via
${\cal H}\colon =({\cal B},{\cal B},S)$.

A right $H$-comodule (algebra) $(A,\tau)$ is an example of a right
${\cal B}$-comodule (algebra) and gives rise to an
${\calH}$-comodule (algebra) via $(A,\tau,\tau)$. The
$H$-coinvariants are obviously 
the same as the $\cal B$-coinvariants and the extension $A^{coH}\subset
A$ is $H$-Galois if and only if it is ${\cal B}$-Galois.
\eex
\bex A weak bialgebra (\cite{Nill,BNSz}) has been shown to determine a
left bialgebroid in (\cite{EtNik,Schau: WHA-qg,Szlachanyi}). Weak
comodule algebras 
(\cite{WeakDoi,CaendeGr00,CaendeGr04}) over a weak bialgebra have been
shown in (\cite{BrzeCaenMil}, Proposition 3.9) to be equivalent to
comodule algebras for the corresponding left bialgebroid. 

Now just the same way as a weak bialgebra determines a left
bialgebroid, it determines also a right bialgebroid, and a weak Hopf
algebra determines a Hopf algebroid (\cite{hgd}, Example 4.8). A weak
comodule algebra for the weak bialgebra is equivalent also to a
comodule algebra for the corresponding right bialgebroid, and a
weak comodule algebra for a weak Hopf algebra is equivalent to a
comodule algebra for the corresponding Hopf algebroid. 

The coinvariants of a weak comodule algebra for a weak bialgebra 
are the same as its coinvariants for the corresponding left or right
bialgebroid. 

A weak bialgebra extension $B \subset A$ is a Galois extension in the sense of  
\cite{CaendeGr04} if and only if $B \subset A$ is a Galois extension by the
corresponding right bialgebroid.  
In the case of a weak Hopf algebra with a bijective antipode,
this is equivalent also to the Galois property of $B \subset A$ as an
extension by the corresponding left bialgebroid.
\eex
\bex The total algebra $H$ of a right bialgebroid
$\calH_R=(H,R,s,t,\gamma,\pi)$ is a right $\calH_R$-comodule algebra
via $\gamma$. We have $H^{co\calH_R}=t(R)$.

The total algebra $H$ of a left bialgebroid
$\calH_L=(H,L,s,t,\gamma,\pi)$ is a right $\calH_L$-comodule algebra
via $\gamma$ and $H^{co\calH_L}=s(L)$.

For a Hopf algebroid $\calH=(\calH_L,\calH_R,S)$, the total algebra is
then a right $\calH$-comodule algebra and the canonical map
$$ {\rm can}_R:{^R H}\ot H_R\to H^R\ot {^RH}\qquad 
h\stac{R\op} h\pri \mapsto h h^{\prime(1)}\stac{R} h^{\prime(2)}$$
is bijective with inverse
$${\rm can}_R\inv:H^R\ot {^RH}\to {^R H}\ot H_R\qquad 
h\stac{R} h\pri \mapsto hS(h\pri\di)\stac{R\op}h\pri\dii.$$
That is, the extension $t_R:R\op\to H$ is $\calH_R$-Galois.
If the antipode $S$ is bijective then also the canonical map
$$ {\rm can}_L:{H^L}\ot {_L H}\to H_L \ot {_L H}\qquad 
h\stac{L} h\pri \mapsto h\di  h^{\prime}\stac{L} h\dii$$
is bijective with inverse
$${\rm can}_L\inv:H_L \ot {_L H}\to {H^L}\ot {_L H}\qquad 
h\stac{L} h\pri \mapsto h^{\prime(2)}\stac{L} S\inv(h^{\prime(1)})h,$$
that is also the extension $s_L:L\to H$ is $\calH_L$-Galois.
\eex
\bex Let $\calH$ be a Hopf algebroid and $A$ a right $\calH_R$-module
algebra. The smash product algebra \cite{KadSzlach} $A\rtimes H$ is
defined as the $k$-module $A^R\ot {^R H}$ with multiplication
$$ (a\rtimes h)(a\pri\rtimes h\pri)\colon = a\pri (a\cdot
h^{\prime(1)})\rtimes h h^{\prime(2)}.$$
With this definition $A\rtimes H$ is an $R$-ring via the homomorphism
$$R\to A\rtimes H\qquad r\mapsto 1_A\rtimes s_R(r) $$
or an $L\op$-ring via the anti-homomorphism
$$L\to A\rtimes H\qquad l\mapsto 1_A\rtimes t_L(l). $$
One can introduce right $\calH_L$- and right $\calH_R$-comodule structures on
$A\rtimes H$ via $\tau_L\colon = A\stac{R} \gamma_L$ and $\tau_R\colon
=A\stac{R} \gamma_R$, respectively. The triple $(A\rtimes H,\tau_L,\tau_R)$ is
a right $\calH$-comodule algebra. We have 
$(A\rtimes H)^{co\calH_R}=\{a\rtimes 1_H\}_{a\in A}$ and the canonical map
\bea
{\rm can}_R:(A\rtimes H)\stac{A}(A\rtimes H)\simeq 
A^R\ot {^RH_{\bf R}}\ot {^{\bf R} H}&\to& (A\rtimes H)\stac{R} H\simeq
A^R\ot  {^RH^{\bf R}}\ot {^{\bf R} H}\nn
a\stac{R} h\stac{\bf R} h\pri &\mapsto& a\stac{R} h\pri h\ui \stac{\bf
R} h\uii\nonumber\eea
is bijective with inverse
$${\rm can}_R\inv: A^R\ot  {^RH^{\bf R}}\ot {^{\bf R} H}\to
A^R\ot {^RH_{\bf R}}\ot {^{\bf R} H}\qquad
a\stac{R} h\stac{\bf R} h\pri \mapsto a\stac{R} h\pri\dii\stac{\bf R}
h S(h\pri\di).$$
This means that the extension $A\subset A\rtimes H$ is
$\calH_R$-Galois. If the antipode of $\calH$ is bijective then it is
also $\calH_L$-Galois.
\eex
\bex In (\cite{Kadison}, Theorem 5.1) Kadison has shown that a depth 2
extension $B\subset A$ of $k$-algebras -- if it is balanced or faithfully
flat -- is a Galois extension for the right bialgebroid, constructed in
\cite{KadSzlach} on the total algebra $(A\stac{B} A)^B$ (the
centralizer of $B$ in the canonical bimodule $A\stac{B} A$).

Recall that if the extension $B\subset A$ is in addition a Frobenius
extension, $(A\stac{B} A)^B$ possesses a Frobenius Hopf algebroid structure
\cite{hgd}. Extending the result of \cite{Kadison} considerably,
B\'alint and Szlach\'anyi have shown in (\cite{BalSzlach}, Theorem
3.7) that an extension $B\subset A$ of $k$-algebras is $\calH$-Galois for
some Frobenius Hopf algebroid $\calH$ if and only if it is a
balanced depth 2 Frobenius extension.
\eex


\section{A Kreimer-Takeuchi type theorem for Hopf algebroids}
\lb{sec:KT}

In this section we investigate $\cal H$-extensions for finitely generated
projective Hopf algebroids $\calH$
with a bijective antipode. We show that for an $\calH$-Galois 
extension $B\subset A$, the algebra $A$ is projective both as a
left and as a right $B$-module and -- under the
additional assumption that $(A\stac{k} A)^{co\calH}\simeq A\stac{k} B$, 
see below -- the surjectivity of the canonical map (\ref{canR})
implies its bijectivity. This is a generalization of
the classical theorem for finitely generated projective Hopf algebras
by Kreimer and Takeuchi (\cite{KreiTak}, Theorem 1.7).

Recently 
Schauenburg and Schneider \cite{SchauSchne} have used new ideas to prove the
Kreimer-Takeuchi theorem and generalizations of it. Their arguments are
formulated in terms of entwining structures \cite{BrzeMaj} over a commutative
ring. In what follows we claim that the line of reasoning in \cite{SchauSchne}
can be applied almost without modification to entwining structures over
non-commutative algebras so to prove a Kreimer-Takeuchi type theorem
for Hopf algebroids.

\bigskip

As we have seen at the beginning of Section \ref{sec:Hopfext}, a right
comodule algebra $A$ for a right $R$-bialgebroid $\calH_R$ determines an
entwining structure (\ref{entw}) over $R$. 
\bl \lb{psibij}
Let $\calH=({\calH}_L,{\calH}_R,S)$ be  a Hopf algebroid with bijective
antipode and $A$ be a right $\calH$-comodule algebra. The map $\psi$ in
(\ref{entw}), corresponding to the right $\calH_R$-comodule algebra
structure of $A$, is an isomorphism. 
\el
\pr The inverse of $\psi$ is constructed using the $\calH_L$-coaction
$a\mapsto a_{\langle 0\rangle} \stac{L} a_{\langle 1\rangle}$ on 
$A$, as
$$\psi\inv:A\stac{R} H\to H\stac{R} A\qquad a\stac{R} h\mapsto
hS\inv(a_{\langle 1\rangle})\stac{R} a_{\langle 0\rangle}.$$
\hfill\qed
\smallskip

Motivated by Lemma \ref{psibij}, we study bijective entwining
structures over $R$. We are going to generalize (\cite{SchauSchne},
Theorem 3.1).
Recall that for any entwined module $M$ and any 
$k$-module $V$, the $k$-module $V\stac{k} M$ is an entwined module via
the second tensor factor. The elements of $V\stac{k} M^{coC}$ form a
subset of $(V\stac{k} M)^{coC}$. We have $V\stac{k} M^{coC}=
(V\stac{k} M)^{coC}$ if, for example, the $k$-module $V$ is flat.
\bp \lb{SchSch}
Let $(A,C,\psi)$ be  a bijective entwining structure over the
algebra $R$, such that the $R$-coring $C$ possesses a grouplike
element $e$. Denote the corresponding right $C$-coaction on $A$ by
$a_{\langle 0\rangle} \stac{R} a_{\langle 1\rangle} \colon
=\psi(e\stac{R} a)$ and denote the subring of its coinvariants by
$B$. Assume that $C$ is flat as a
left (right) $R$-module and projective as a right (left) $C$-comodule.

{\em (1)} Suppose that $(A\stac{k} A)^{coC}=A\stac{k} B$ and the 
canonical map 
\be {\rm can}: A\stac{B} A\to A\stac{R} C\qquad a\stac{B} a\pri\mapsto
aa\pri_{\langle 0\rangle}\stac{R}a\pri_{\langle 1\rangle}\lb{can}\ee
is surjective. Under these assumptions the canonical map (\ref{can}) is
bijective. 

{\em (2)} If the canonical map (\ref{can}) is bijective then $A$
is projective as a right (left) $B$-module.
\ep
\pr The proof is actually the same as the proof of (\cite{SchauSchne},
Theorem 3.1), so we present only a sketchy proof here.

{\underline{(1)}} Let us use the assumption that $C$ is projective as
a right $C$-comodule.
By the bijectivity of $\psi$ and the adjunction (\ref{adjA}),
for any 
entwined module $M$ we have $\Hom^C_A(A\stac{R} C,M)\simeq
\Hom^C(C,M)$. 
The forgetful functor ${\mathcal M}^C_A \to {\mathcal M}^C$ is a right
adjoint, hence it preserves monomorphisms. A morphism in any of the (comodule)
categories ${\mathcal M}^C_A$ and ${\mathcal M}^C$ is an epimorphism if and
only if it is a surjective map, hence the forgetful functor ${\mathcal M}^C_A
\to {\mathcal M}^C$ preserves also epimorphisms. 
Therefore, in light of (\cite{BrzeWis}, 18.20(2)), by 
flatness of the left $R$-module $C$, projectivity of the right
$C$-comodule $C$ implies projectivity of the entwined module
$A\stac{R} C$. Hence the surjective map
\begin{equation}\label{eq:lifted}
A\stac{k} A\to A\stac{R} C\qquad a\stac{k} a\pri \mapsto 
aa\pri_{\langle 0\rangle}\stac{R}a\pri_{\langle 1\rangle}
\end{equation}
is a split epimorphism of entwined modules. This means that $A\stac{R} C$ is a 
direct summand of $A\stac{k} A$.

Notice that $(A\stac{R} C)^{coC}\simeq A$ via the isomorphism
$$\alpha: A\to (A\stac{R} C)^{coC}\qquad a\mapsto a \stac{R} e,$$
hence the canonical map (\ref{can}) is related to the unit of the
adjunction (\ref{adj}) as ${\rm can}=\mu_{A\stac{R} C}\ci
(\alpha\stac{B} A)$. This means that ${\rm can}$ is bijective provided
$\mu_{A\stac{k}A}$ is bijective.
Tensoring a $k$-free resolution 
$$\dots P_n \stackrel{\partial_n}{\longrightarrow} P_{n-1} \longrightarrow
\dots \longrightarrow P_1 \stackrel{\partial_1}{\longrightarrow} P_0
\stackrel{\partial_0}{\longrightarrow} A\longrightarrow 0$$ 
of $A$ with $A$ over $k$, we can write $A\stac{k} A$ as
the cokernel of the morphism  $\partial_1\stac{k} A: P_1\stac{k} A\to 
P_0\stac{k} A$ between entwined modules that are direct sums of copies of
$A$. Functoriality of $\mu$ implies that 
$$\mu_{A\stac{k} A}\ci [(\partial_0\stac{k} B)\stac{B} A]=(\partial_0\stac{k}
A) \ci (P_0\stac{k} \mu_A).$$
Since $\mu_A$ is an isomorphism, we can use the universality of the
cokernel to define the inverse of $\mu_{A\stac{k} A}$ as the unique map
$\phi:A\stac{k} A\to (A\stac{k} A)^{coC}\stac{B} A=(A\stac{k}
B)\stac{B} A$ which satisfies
$$ \phi\ci (\partial_0\stac{k} A) = [(\partial_0\stac{k} B)\stac{B} A]\ci
(P_0\stac{k} \mu_A\inv). $$
This proves the bijectivity of $\mu_{A\stac{k} A}$, hence of the canonical map 
(\ref{can}).

If $C$ is projective as a left $C$-comodule and flat as a right
$R$-module then $\psi$ has to be
replaced with its inverse in the above line of reasoning.

{\underline{(2)}}
Projectivity of the right $B$-module $A$ is proven from the
projectivity of the right $C$-comodule $C$ as follows. By
bijectivity of the canonical map (\ref{can}) and using the adjunction
(\ref{adj}), 
for any entwined module $M$ we have $\Hom^C(C,M)\simeq \Hom_B(A,
M^{coC})$. Hence the functor $\Hom_B(A,(\h)^{coC}):\M^C_A\to \M_k$ is
exact.

Let $f:\oplus_I B\to A$ be an epimorphism in $\M_B$, for some index set
$I$. Then $f\stac{B} A$ is an epimorphism in $\M^C_A$ which is mapped
by the functor $\Hom_B(A,(\h)^{coC})$ to the epimorphism $\Hom_B(A,f)$
in $\M_k$. 

In order to prove projectivity of the left $B$-module $A$ from
projectivity of the left $C$-comodule $C$, replace $\psi$ by its
inverse in the above arguments.
\hfill\qed
\smallskip

As it has been explained to us by Tomasz Brzezi\'nski, there exists an
alternative, more general argument proving bijectivity of the
canonical map (\ref{can}) from split surjectivity of (\ref{eq:lifted}) in
$\M^C_A$. Namely, application of (\cite{BTW}, Theorem 2.1) to the $A$-coring
$A\stac{R} C$, corresponding to the $R$-entwining structure
$(A,C,\psi)$ with grouplike element $e\in C$, and its comodule $M=A$,
implies the claim. Indeed, in this case condition b) of (\cite{BTW}, Theorem
2.1) reduces to $(A\stac{k} A)^{co C}\simeq A\stac{k} A^{coC}$.

We have formulated Proposition \ref{SchSch} (1) in terms of the
assumption $(A\stac{k} A)^{coC}=A\stac{k} A^{coC}$. It has the
advantage that in certain situations (e.g. if $k$ is a field) it
obviously holds true. We do not know, however, whether 
it is also a necessary condition for the claim of Proposition
\ref{SchSch} (1). On the other hand, notice that if the
canonical map (\ref{can}) is an isomorphism in ${\M^C_A}$ then
\be (A\stac{B} A)^{coC}\simeq (A\stac{R} C)^{coC}=A,\lb{Rcond}\ee
where $B$ stands for $A^{coC}$, as before.
Application of (\cite{Wisbauer}, Proposition 5.1) to the $A$-coring
$A\stac{R} C$, corresponding to the $R$-entwining structure
$(A,C,\psi)$ with grouplike element $e\in C$, and its comodule $M=A$,
implies that the bijectivity of the canonical map (\ref{can}) follows
from split surjectivity of (\ref{eq:lifted}) also under the (sufficient and
necessary) condition (\ref{Rcond}).
\bigskip

Let us turn to the application of Proposition \ref{SchSch} to Hopf
algebroid extensions. Let $\calH$ be a finitely generated projective
Hopf algebroid with a bijective antipode. 
It follows from the Fundamental Theorem for Hopf algebroids (\cite{integral},
Theorem 4.2, see also the {\em Corrigendum}) 
and the existence of a Hopf module structure on
$H^*=\Hom_R(H,R)$ with coinvariants $\il(H^*)$
(\cite{integral}, Proposition 4.4), that for such a Hopf algebroid the map
\be \lb{alphaL}
\alpha_L: {^L H}\ot \il(H^*)^L\ \to\ H^*\qquad h\stac{L\op} \lambda^*\ \mapsto
\lambda^*\lu S(h)\equiv \lambda^*[S(h)\h]\ee
is an isomorphism of Hopf modules, hence in particular of left $H$-modules. 
(The left $H$-module structure on ${^L H}\ot \il(H^*)^L$ is given by left
multiplication in the first tensor factor, and on $H^*$ by $h\cdot \phi^*\colon
=\phi^*\lu S(h)\equiv \phi^*[S(h)\h]$.)
This implies that the element
$\alpha_L\inv(\pi_R)$ is an invariant of the left $H$-module ${^L
H}\ot \il(H^*)^L$. The elements $\{x_i\}\subset H$ and
$\{\lambda_i^*\}\subset \il(H^*)$ satisfying $\sum_i x_i\stac{L\op}
\lambda^*_i=\alpha_L\inv (\pi_R)$ can be
used to construct dual bases for the left ${^*H}$-module on $H$
defined as 
${^*\phi}\rd h\colon =h\ui \ s_R\!\ci\! {^*\phi}(h\uii)$. As a matter of
fact, for $\lambda^*\in \il(H^*)$ we have $\lambda^*\ci S\in \il({^*H})$ 
(\cite{integral}, Scholium 2.10), and ${^* H}$ is a right $H$ module
via ${^*\phi}\ld h={^*\phi}(h\h)$. We leave it
to the reader to check that the sets $\{x_i\}\subset H$ and 
$\{\lambda^*_i\ci S\ld S\inv(\h)\}\subset {_{^*H}\Hom}(H,{^*H})$ are
dual bases, showing that $H$ is a finitely generated and projective
${^* H}$-module.  

Since the antipode was assumed to be bijective, we can apply the same
argument to the co-opposite Hopf algebroid $\calH_{cop}$ to conclude
on the finitely generated projectivity of the left $H^*$-module on $H$
defined as $\phi^* \ru h\colon=h\uii \ t_R\!\ci\! \phi^*(h\ui)$. 

Furthermore, projectivity of $H$ as a left ${^*H}$ module implies that it is
projective as a right $\calH_R$-comodule and projectivity of $H$ as a left
$H^*$-module implies that it is projective as a left
$\calH_R$-comodule. 
These observations allow for the application of Proposition
\ref{SchSch} to the entwining structure (\ref{entw}) -- which is
bijective by Lemma \ref{psibij}. 
\bc\lb{KT}
Let $\cal H$ be a finitely generated projective Hopf algebroid with a
bijective antipode and let $B\subset A$ be an $\calH$-extension. 

{\em (1)} Suppose that $(A\stac{k} A)^{co\calH}=A\stac{k} B$ (e.g. $A$ is
$k$-flat). 
If the canonical map 
$${\rm can}_R: A\stac{B} A \to A\stac{R} H \qquad a\stac{B} a\pri
\mapsto a a^{\prime \langle 0 \rangle} \stac{R} a^{\prime \langle 1
\rangle}$$ 
is surjective then the extension $B\subset A$ is $\calH_R$-Galois
(equivalently, $\calH_L$-Galois). 

{\em (2)} If the extension $B\subset A$ is $\calH$-Galois then $A$ is
projective both as a left and as a right $B$-module.
\ec


\section{Morita theory for Hopf algebroid extensions}
\lb{sec:Morita}

As we have seen at the beginning of Section \ref{sec:Hopfext}, an
$\calH_R$-extension $B\subset A$, for a right $R$ -bialgebroid  $\calH_R$,
determines an $A$-coring structure on $A\stac{R} H$. One 
can apply the Morita theory for corings, developed in
\cite{CaenVerWang}, to this coring. In particular, if $\calH_R$ is
finitely generated and projective as a left $R$-module, then Theorem
\ref{CVW} (1) can be used to obtain criteria which are equivalent to
the $\calH_R$-Galois 
property of the extension $B\subset A$ and the projectivity of the
left $B$-module $A$. 
Analogously, one can apply 
Theorem \ref{CVW} (1) to obtain criteria which are equivalent to the Galois
property of an $\calH_L$-extension $B\subset A$ for a finitely
generated projective left bialgebroid
$\calH_L$ and the projectivity of the right $B$-module $A$. 

Applying the results of the previous sections, we prove that if
$\calH=(\calH_L,\calH_R,S)$  
is a finitely generated projective Hopf algebroid with a bijective
antipode and 
$A$ is a right $\calH$-comodule algebra, with $\calH_R$, equivalently,
$\calH_L$-coinvariant subalgebra $B$,
then the equivalent conditions, derived from Theorem \ref{CVW} (1)
for $B\subset A$ as an $\calH_R$-extension, on one hand, and as an
$\calH_L$-extension, on the other hand,
are equivalent also to each other.

\bigskip

Let $\calH_R=(H,R,s,t,\gamma,\pi)$ be a right bialgebroid such that $H$ is
finitely generated and projective as a left $R$-module and let $A$ be
a right $\calH_R$-comodule algebra.
As a first step, let us identify the Morita context $(A^{^*C}, {^*C},A,
{^*C}^{^*C},\nu,\mu)$, associated to the $A$-coring $C=A\stac{R} H$
(as it is explained in Section \ref{ss:Morcor}).

Recall that $A$ -- being a right $\calH_R$-comodule -- is also a left
${^*H}$-module via 
${^*\phi}\cdot a=a^{\langle 0\rangle}\cdot {^*\phi}(a^{\langle 1\rangle})$.
Since the module ${^R H}$ is finitely generated and
projective, ${^*H}$ possesses a left bialgebroid structure.
Let us introduce the smash product algebra ${^*H}\ltimes A$ as the 
$k$-module ${^*H} \stac{R} A$ (where the right $R$-module structure of
${^*H}$ is given by $({^*\phi}\cdot r)(h)={^*\phi}(h)r$) with
multiplication 
\be ({^* \phi}\ltimes a)({^*\psi}\ltimes a\pri)\colon = {^*\psi}\di
{^*\phi}\ltimes ({^*\psi}\dii\cdot a)a\pri.
\lb{smash}\ee
With this definition ${^*H}\ltimes A$ is an $A$-ring via
the homomorphism
$$i_A:A\to {^* H}\ltimes A\qquad a\mapsto \pi\ltimes a.$$
Define a right-right {\em relative $(A,\calH_R)$ -module} to be
an entwined module for the $R$ -entwining structure (\ref{entw}),
i.e. a right comodule for the $A$-coring $A\stac{R} H$. This means a
right $A$-module (and hence in particular a right $R$-module) and a
right $\calH_R$-comodule $M$ such that the compatibility condition
$$(m\cdot a)^{\langle 0\rangle}\stac{R} (m\cdot a)^{\langle 1\rangle}=
m^{\langle 0\rangle}\cdot a^{\langle 0\rangle}\stac{R}
m^{\langle 1\rangle}a^{\langle 1\rangle}$$
holds true for any $m\in M$ and $a\in A$. 
Clearly, $A$ is itself a relative $(A,\calH_R)$-module. 
It is straightforward to check that the category $\M^H_A$ of relative
$(A,\calH_R)$-modules is isomorphic to the category 
$\M_{{^* H}\ltimes A}$ of right ${^* H}\ltimes A$-modules and the ${^*
H}\ltimes A$-invariants are the same as the
$\calH_R$-coinvariants. This fact can be easily understood in view of
\bl \lb{MorHopf}
Let $\calH_R$ be a right $R$-bialgebroid such that $H$ is
finitely generated and projective as a left $R$-module and let $A$ be
a right $\calH_R$-comodule algebra.
Then the left $A$-dual algebra of the $A$-coring $A\stac{R} H$ is
isomorphic to the smash product algebra ${^* H}\ltimes A$. 
\el
\pr The required algebra isomorphism is constructed as
$${^*H}\ltimes A\to {_A\Hom}(A\stac{R} H, A)\qquad 
{^*\phi}\ltimes a\mapsto (\ a\pri \stac{R} h\mapsto a\pri
({^*\phi}(h)\cdot a)\ ).$$
\hfill\qed
\smallskip

The Morita context, associated to the $\calH_R$-extension $B\subset A$,
is then $(B,{^*H}\ltimes A, A, ({^*H}\ltimes A)^{co\calH_R}$ 
{\mbox{$\nu_R,\mu_R)$ }}
with connecting maps
\bea 
&\nu_R: A\stac{{^*H}\ltimes A} ({^* H}\ltimes A)^{co\calH_R}\to B
\qquad & a\pri \stac{{^*H}\ltimes A}(\sum_i {^*\phi}_i\ltimes a_i)
\mapsto \sum_i ({^*\phi}_i\cdot a\pri) a_i\\
&\mu_R: ({^* H}\ltimes A)^{co\calH_R} \stac{B} A\to 
{^* H}\ltimes A  \qquad 
&(\sum_i {^*\phi}_i\ltimes a_i)\stac{B} a\pri \mapsto \sum_i
{^*\phi_i}\ltimes a_i a\pri.\lb{muR}
\eea
Since $H$ is finitely generated and projective as a left $R$-module, so
is the left $A$-module $A\stac{R} H$. Hence part (1) in Theorem
\ref{CVW} implies
\bp \lb{right}
Let $\calH_R$ be a right $R$-bialgebroid such that $H$ is
finitely generated and projective as a left $R$-module and let $B\subset A$ be
an $\calH_R$-extension. The following assertions are equivalent.

{\em (a)} The map $\mu_R$ in (\ref{muR}) is surjective (and, a
fortiori, bijective).

{\em (b)} The functor $(\h)^{co\calH_R}:\M^H_A\to \M_B $ is fully
faithful.

{\em (c)} $A$ is a right ${^* H}\ltimes A$-generator.

{\em (d)} $A$ is projective as a left $B$-module and the map
\be {^* H}\ltimes A\to {_B\End}(A)\qquad {^*\phi}\ltimes a\mapsto (\
a\pri\mapsto ({^*\phi}\cdot a\pri)a\ )\lb{*can}\ee
is an algebra anti-isomorphism.

{\em (e)} $A$ is projective as a left $B$-module and the extension
$B\subset A$ is $\calH_R$-Galois.
\ep
The arguments leading to Theorem \ref{right} can be repeated by replacing 
the right bialgebroid $\calH_R$ with a left $L$-bialgebroid $\calH_L$ such
that $H$ is finitely generated and projective as a left $L$-module. Indeed, 
in this case the left $L$-dual, ${_*H}$, possesses a right bialgebroid
structure 
and $A$ is a right ${_*H}$-module via $a\cdot {_*\phi}=
a_{\langle 0\rangle}\cdot {_*\phi}(a_{\langle 1\rangle})$.
The right $A$-dual of the $A$-coring $A\stac{L} H$ is isomorphic
to the smash product algebra ${_* H}\ltimes A$, which is
defined as the $k$-module ${_*H}\stac{L} A$, (where the right $L$-module
structure on ${_* H}$ is given by $({_*\phi}\cdot l)(h)={_*\phi}(h)l$), with
multiplication  
$$ ({_*\phi}\ltimes a)({_*\psi}\ltimes a\pri)\colon = 
{_*\psi}{_*\phi}\ui \ltimes a (a\pri\cdot {_*\phi}\uii).$$
Since $A$ is an $L\op$-ring, a left $A$-module is in particular a right
$L$-module and we have the isomorphism $(A\stac{L} H)\stac{A} M\simeq
M\stac{L} H$ for any left $A$-module $M$. A left comodule for the $A$-coring
$A\stac{L} H$ is then equivalent to a left $A$-module (hence in particular a
right $L$-module) and a right $\calH_L$-comodule $M$ such that the
compatibility condition
$$ (a\cdot m)_{\langle 0\rangle}\stac{L} (a\cdot m)_{\langle 1\rangle}=
a_{\langle 0\rangle}\cdot m_{\langle 0\rangle}\stac{L} 
a_{\langle 1\rangle} m_{\langle 1\rangle}$$
holds true for $m\in M$ and $a\in A$. Such modules are called left-right
relative $(A,\calH_L)$-modules 
and their category is denoted by  ${_A\M^H}$. It follows that the category
${_A\M^H}$ is isomorphic also to
${_{{_*H}\ltimes A}\M}$, the category of left ${_*H}\ltimes A$-modules. 
The Morita context associated to the $\calH_L$-extension $B\subset A$ is $(B,
{_* H}\ltimes A,A, ({_* H}\ltimes A)^{co\calH_L},\nu_L,\mu_L)$ with connecting
maps 
\bea
&\nu_L:  ({_* H}\ltimes A)^{co\calH_L}\stac{{_* H}\ltimes A} A\to B\qquad
&(\sum_i {_*\phi}_i\ltimes a_i)\stac{{_* H}\ltimes A} a\pri \mapsto \sum_i
a_i(a\pri\cdot {_*\phi}_i)\\
&\mu_L:A\stac{B} ({_* H}\ltimes A)^{co\calH_L} \to {_* H}\ltimes A\qquad 
&a\pri \stac{B} (\sum_i {_*\phi}_i\ltimes a_i)\mapsto \sum_i {_*\phi}_i\ltimes
a\pri a_i\lb{muL}.
\eea
Part (1) of Theorem \ref{CVW} implies
\bp\lb{left}
Let $\calH_L$ be a left $L$-bialgebroid such that $H$ is
finitely generated and projective as a left $L$-module and let $B\subset A$ be
an $\calH_L$-extension. The following assertions are equivalent.

{\em (a)} The map $\mu_L$ in (\ref{muL}) is surjective (and, a
fortiori, bijective).

{\em (b)} The functor $(\h)^{co\calH_L}:{_A\M^H}\to {_B \M} $ is fully
faithful.

{\em (c)} $A$ is a left ${_* H}\ltimes A$-generator.

{\em (d)} $A$ is projective as a right $B$-module and the map
\be {_* H}\ltimes A\to {\End_B}(A)\qquad {_*\phi}\ltimes a\mapsto (\
a\pri\mapsto a (a\pri\cdot {_*\phi})\ )\lb{can*}\ee
is an algebra isomorphism.

{\em (e)} $A$ is projective as a right $B$-module and the extension
$B\subset A$ is $\calH_L$-Galois.
\ep
Combining Proposition \ref{right}, \ref{left}, Corollary \ref{KT} and Lemma
\ref{LReq} we can state our main result.
\bt \lb{CFM}
Let $\calH=(\calH_L,\calH_R,S)$ be a finitely generated projective Hopf
algebroid with a bijective antipode and let $B\subset A$ be an
$\calH$-extension. The following assertions are equivalent.

{\em (a)} The extension $B\subset A$ is $\calH_R$ -Galois.

{\em (b)} $A$ is projective as a left $B$-module and the map (\ref{*can})
is an algebra anti-isomorphism.

{\em (c)} $A$ is a right ${^* H}\ltimes A$-generator.

{\em (d)} The functor $(\h)^{co\calH_R}:\M^H_A\to \M_B $ is fully
faithful.

{\em (e)} The map $\mu_R$ in (\ref{muR}) is surjective (and, a
fortiori, bijective).

{\em (a$\pri$)} The extension $B\subset A$ is $\calH_L$ -Galois.

{\em (b$\pri$)} $A$ is projective as a right $B$-module and the map
(\ref{can*}) is an algebra isomorphism.

{\em (c$\pri$)} $A$ is a left ${_* H}\ltimes A$-generator.

{\em (d$\pri$)} The functor $(\h)^{co\calH_L}:{_A\M^H}\to {_B \M} $ is fully
faithful.

{\em (e$\pri$)} The map $\mu_L$ in (\ref{muL}) is surjective (and, a
fortiori, bijective).
\et
\pr It follows from part (2) of Corollary \ref{KT} that (a) is equivalent to
any of the assertions
\bea 
&&A\ \textrm{is\ projective\ as\ a\ left\ }B\ \textrm{module\ and\ the\
  extension\ }B\subset A\ \textrm{is\ }\calH_R\ \textrm{-Galois}.\lb{*}\\
&&A\ \textrm{is\ projective\ as\ a\ right\ }B\ \textrm{module\ and\ the\
  extension\ }B\subset A\ \textrm{is\ }\calH_R\ \textrm{-Galois}.\lb{**}
\eea
By Lemma \ref{LReq} assertions (a) and (a$\pri$) are equivalent and
(\ref{**}) is equivalent to
\be
A\ \textrm{is\ projective\ as\ a\ right\ }B\ \textrm{module\ and\ the\
  extension\ }B\subset A\ \textrm{is\ }\calH_L\ \textrm{-Galois}.\lb{***}\ee
The rest of the proof follows from Propositions \ref{right} and \ref{left}.
\hfill\qed
\smallskip

Let us mention that in (\cite{Caenepeel}, Theorem 4.7) a stronger version of
Theorem \ref{CVW} (1) has been proven. Its application to bialgebroid
extensions implies
\bp\lb{sright}
Let $\calH_R$ be a right $R$- bialgebroid such that $H$ is
finitely generated and projective as a left $R$-module and let $B\subset A$ be
an $\calH_R$-extension. The following assertions are equivalent.

{\em (a)} The Morita context $(B,{^* H}\ltimes A,A,({^* H}\ltimes
A)^{co\calH_R},\nu_R,\mu_R)$ is strict. 

{\em (b)} The functor $(\h)^{co\calH_R}:\M^H_A\to \M_B $ is an equivalence
with inverse $\h\stac{B} A: \M_B \to \M^H_A$.

{\em (c)}  The map (\ref{*can}) is an algebra anti-isomorphism and 
$A$ is a left $B$-progenerator.

{\em (d)} $A$ is faithfully flat as a left $B$-module and the extension
$B\subset A$ is $\calH_R$-Galois.
\ep
Also the analogue of Proposition \ref{sright} for left bialgebroid extensions
can 
be proven. We can not prove, however, that for an $\calH$-extension,
in the case of any finitely generated projective Hopf algebroid $\calH$
with a bijective antipode, the equivalent
conditions in Proposition \ref{sright} and their counterparts on the left
bialgebroid of $\calH$ are equivalent also to each other (as it was
seen to be   
the case with Proposition \ref{right} and \ref{left}).
\footnote[3]{This question is answered in \cite{ABM}, see Proposition 16 in
  the {\em Corrigendum} or Proposition 4.4 in the arXiv version.}
\bigskip

The Morita context $(A^{^*C}, {^*C},A,{^*C}^{^*C},\nu,\mu)$, associated in
\cite{CaenVerWang} to an $A$-coring $C$ with a grouplike element, is 
a generalization of the Morita context, associated to a bialgebra extension
in \cite{Doi}. In the case of a finite dimensional Hopf algebra over a field
(or 
a Frobenius Hopf algebra over a commutative ring) another Morita context has
been associated to a Hopf algebra extension in \cite{CohFishMont,CohFish}. The
relation of the two Morita contexts is of the type described in part (2) of
Theorem \ref{CVW}.
In order to see what is the analogue of the Morita context of Cohen, Fishman
and Montgomery in the case of Hopf algebroids, in the rest of the section we
assume that $\calH$ is a {\em Frobenius} Hopf algebroid.
\bl \lb{frob}
Let $\calH$ be a  Frobenius Hopf algebroid and $A$ be a left $\calH_L$-module
algebra. Consider the smash product algebra $H\ltimes A$, which is the
$k$-module $H\stac{L} A$ with multiplication 
$$ (h\ltimes a)(g\ltimes a\pri)\colon = g\di h\ltimes (g\dii\cdot a)a\pri.$$
The extension 
$$ i:A\to H\ltimes A\qquad a\mapsto 1_H\ltimes a $$
is a Frobenius extension.
\el
\pr Recall (from Section \ref{ss:hgd}) that a Frobenius Hopf algebroid
possesses non-degenerate left integrals. Let us fix such an integral $\ell$
and denote by $\rho_*$ the unique element in $H_*$, for which 
$\ell\lu\rho_*\equiv s_L\ci \rho_*(\ell\di)\ell\dii=1_H$. 
A Frobenius functional $\Phi:H\ltimes A\to A$ is
given by $h\ltimes a\mapsto \rho_*(h)\cdot a$. A Hopf algebroid calculation
shows that it is an $A$-$A$ bimodule map and possesses a dual basis
$(S(\ell\uii)\ltimes 1_A)\stac{A}(\ell\ui\ltimes 1_A)$.
\hfill\qed
\smallskip

Recall that for a Frobenius Hopf algebroid $\calH$ also the left bialgebroid
${^*H}$ 
possesses a Frobenius Hopf algebroid structure.
Applying Lemma \ref{frob} together with part (2) of Theorem \ref{CVW} 
we conclude that the Morita context $(A^{co\calH},{^*H}\ltimes A,
A,({^*H}\ltimes A)^{co\calH},\nu_R,\mu_R)$, associated to the right
$\calH_R$ -comodule algebra structure of a right $\calH$-comodule algebra
$A$ for the Frobenius Hopf algebroid $\calH$, is 
equivalent to the Morita context $(A^{co\calH},{^*H}\ltimes A,
A,A,\nu\pri,\mu\pri)$ with connecting maps 
\bea 
&\nu\pri:A\stac{{^*H}\ltimes A} A\to A^{co\calH}\qquad 
&a\stac{{^*H}\ltimes A} a\pri \mapsto {^*\lambda}\cdot (aa\pri)\lb{nupri}\\
&\mu\pri: A\stac{A^{co\calH}} A\to {^*H}\ltimes A\qquad
&a\stac{A^{co\calH}} a\pri \mapsto (\pi_R\ltimes a)({^*\lambda} \ltimes 1_A)
(\pi_R\ltimes a\pri),
\lb{mupri}\eea
where ${^*\lambda}$ is a non-degenerate left integral in ${^*H}$.
\bc If we add to the conditions of Theorem \ref{CFM} the requirement that
$\calH$ be a Frobenius Hopf algebroid then we can add to the equivalent
assertions (a)-(e$\pri$) also

(f) For any non-degenerate left integral ${^*\lambda}$ in the Frobenius Hopf
  algebroid ${^*H}$ the map
$$A\stac{B} A\to {^*H}\ltimes A\qquad a\stac{B} a\pri\mapsto (\pi_R\ltimes
  a)({^*\lambda}\ltimes 1_A)(\pi_R\ltimes a\pri)$$
is surjective (and, a fortiori, bijective).
\ec
By (\cite{hgd}, Lemma 5.14), for any non-degenerate left integral $\ell$ in a
Hopf algebroid $\calH$, $\rho_*\colon = \ell_L\inv(1_H)\in {H_*}$,
and any element $h$ of $H$, we have $\ell\di h\stac{L} \ell\dii=\ell\di
\stac{L} \ell\dii t_L\ci {\rho_*}(\ell h\ui)S(h\uii)$. This implies that the
image of the map (\ref{mupri})
is an ideal in ${^*H}\ltimes A$. Hence -- just as it has been proven for finite
dimensional Hopf algebras in (\cite{CohFishMont}, Corollary 1.3) -- we see
that it is true also for $\calH$-extensions $B\subset A$ for a
Frobenius Hopf algebroid $\calH$ that if the 
$k$-algebra ${^*H}\ltimes A$ is simple then the extension $B \subset A$ is
$\calH$-Galois. 


\end{document}